# FAIRE DES MATHÉMATIQUES À L'ÉCOLE MATERNELLE : À QUELLES CONDITIONS ?

**Magali HERSANT**[1]

INSPÉ de l'Académie de Nantes, CREN, Nantes Université

**Résumé.** Dans cet article, nous identifions des conditions didactiques à l'apprentissage des mathématiques à l'école maternelle. Nous nous appuyons pour cela sur le cadre de la théorie des situations didactiques (Brousseau, 1998) et la notion de situation-problème (Douady, 1984). Nous explicitons d'abord ce qui constitue pour nous les enjeux de l'enseignement des mathématiques à l'école maternelle puis, à partir d'exemples, nous mettons en évidence des conditions relatives aux enjeux de l'activité des élèves, aux caractéristiques des situations proposées aux élèves et aux interventions de l'enseignant.

**Mots-clés.** Maternelle, situation-problème, apprentissage par adaptation, institutionnalisation, fiche.

## Introduction

L'école maternelle constitue le lieu de la première rencontre des élèves avec l'école et, en particulier avec les apprentissages mathématiques. Des travaux de sociologues et de didacticiens (Bautier, 2006 ; Laparra & Margolinas, 2016) montrent que cette école est parfois aussi déjà un lieu de création des premières inégalités scolaires, y compris dans les apprentissages mathématiques. Notre expérience de formatrice et les nombreuses observations que nous avons effectuées dans des classes nous conduisent à penser que les enseignants peinent parfois à identifier les conditions didactiques génératrices d'apprentissages mathématiques à l'école maternelle, ce qui peut ensuite être source de malentendus scolaires pour les élèves (Bautier & Rayou, 2009) et d'inégalités potentielles face aux apprentissages (Hersant, 2020). Par ailleurs, les ressources disponibles n'aident pas toujours à cet effet dans la mesure où, souvent, peu de commentaires didactiques y sont proposés. À travers cet article qui s'appuie sur l'analyse de situations observées en classe, nous souhaitons mettre en évidence certaines de ces conditions didactiques afin de permettre des apprentissages pour tous les élèves.

Précisons d'emblée que les conditions que nous identifions découlent principalement de deux hypothèses sur lesquelles s'accordent les travaux dans le champ de la didactique des mathématiques et, en particulier, le cadre de la théorie des situations didactiques de Brousseau (1997, 1998). La première correspond à une hypothèse d'apprentissage par adaptation qui prend appui sur les travaux de Piaget : l'élève construit des connaissances en se confrontant à des situations problématiques qui lui renvoient de l'information sur ses actions, ce qui lui permet d'adapter ses connaissances et ses actions pour les situations ultérieures. La seconde hypothèse renvoie aux travaux de Bruner et Vygotski et rend compte de la dimension sociale des apprentissages : l'action permet l'émergence de connaissances, personnelles et souvent implicites, mais les savoirs sont de l'ordre de la culture ; il est donc nécessaire de permettre à ces connaissances de prendre un statut de savoir, dépersonnalisé et partagé dans la classe, à travers des formulations qui contribuent au processus d'institutionnalisation (Brousseau, 1998).

---

[1] magali.hersant@univ-nantes.fr

## 1. Quels enjeux pour l'enseignement des mathématiques à la maternelle ?

### 1.1. La question des savoirs

L'identification des savoirs qui doivent faire l'objet d'apprentissages à l'école maternelle apparaît comme une condition indispensable, d'une part, pour décider de la pertinence de proposer une situation au regard des apprentissages déjà réalisés par les élèves et, d'autre part, pour permettre un processus d'institutionnalisation. Mais ce point est délicat. Laparra et Margolinas (2016) relèvent en effet une certaine transparence des savoirs à la maternelle, c'est-à-dire que les choix des situations et des activités ne sont pas déterminés par les savoirs que portent ces activités et situations mais par des objectifs plus larges, comme par exemple, l'importance de faire manipuler les élèves. Selon nous, deux raisons (au moins) peuvent expliquer ce phénomène. D'abord les savoirs en jeu à l'école maternelle sont tellement incorporés et familiers pour un adulte qu'on peut arriver à les penser naturels et innés. Ensuite, ces savoirs ne sont pas, le plus souvent, explicités dans les manuels et, dans les programmes (MEN, 2021), ils restent masqués par la formulation des objectifs en termes de compétences. On trouve en effet dans ces ressources des objectifs plus larges qui, de notre point de vue, renferment des savoirs fins mais ne les explicitent pas. Ainsi, les programmes identifient des « attendus » sur l'utilisation du nombre qui correspondent à des objectifs en termes de types de tâches mais ne mentionnent pas les savoirs qui interviennent pour leur réalisation. Par exemple, *« réaliser une collection dont le cardinal est compris entre 1 et 10 »* (MEN (2021) implique des savoirs fins très différents selon que le cardinal de la collection est donné avec le mot-nombre énoncé oralement, avec l'écriture chiffrée du nombre ou encore à l'aide d'une décomposition. Les ouvrages annoncent souvent les objectifs associés aux situations qu'ils proposent. On trouve par exemple des objectifs comme *« prendre conscience que les collections équipotentes à une collection donnée sont équipotentes entre elles »* (Valentin, 2004, p. 51) ou *« réaliser une collection dont le cardinal est donné »* (Duprey *et al.*, 2010, p. 36). Mais, là encore, ces objectifs ne disent rien des savoirs fins comme l'énumération (Briand *et al.*, 1999), les procédures de dénombrement ou la connaissance de faits numériques ou de constellations, que l'on souhaite que les élèves construisent ou mobilisent dans une situation ; ils n'aident pas non plus à anticiper ce que l'on peut dire aux élèves pour les aider en cours d'activité, ni à prévoir ce que l'on peut dire, de façon réaliste, à un élève de maternelle sur ce qu'il est en train d'apprendre et qui pourra lui être utile pour les situations ultérieures. On se heurte donc à une double difficulté : les termes savants utilisés pour communiquer sur des objectifs d'apprentissage à long terme ou des catégories de problèmes mathématiques (comparer des collections, constituer des collections équipotentes…) ne permettent pas aux enseignants de repérer les savoirs à apprendre ; par ailleurs ils ne sont pas adaptés pour communiquer avec les élèves, en situation, sur les connaissances en jeu.

Une solution pour aider à identifier ces savoirs fins peut être de chercher ce qu'on peut dire pour aider les élèves à réaliser une tâche donnée. On accède ainsi à des formulations de ces savoirs fins, plus ou moins contextualisées. Par exemple, dans une situation de dénombrement d'une collection d'objets constituée de briques de différentes couleurs, dire que « il y a deux briques rouges et deux briques bleues, en tout il y a quatre briques » ou « deux ici et encore deux là, c'est quatre » sont des éléments de savoirs qui contribuent à l'apprentissage du nombre comme cardinal d'une collection, à la connaissance d'un fait numérique (2+2=4) et sont constitutifs d'une procédure de dénombrement. Des formulations proches peuvent être utilisées dans une situation de comparaison de hauteurs de tours ; elles contribueront alors à l'apprentissage de procédures de comparaison de collections. On le voit, il n'est nullement besoin d'énoncer des objectifs généraux qui ne parleraient pas aux élèves.

Évidemment, creuser la question des savoirs amène à se demander, comme nous l'a formulé une enseignante de PS : « Mais, au fond du fond, qu'est-ce qu'il y a à apprendre sur le nombre à l'école maternelle ? ». Autrement dit, effectivement, il faut bien réussir à identifier des apprentissages « plus gros », comme *« savoir constituer une collection équipotente à une collection donnée lorsque les deux collections sont éloignées »* (*ibid.*), pour se donner des lignes directrices et penser une progression qui n'élude aucun point essentiel.

Concernant les aspects numériques, une part importante des apprentissages à réaliser à l'école maternelle réside dans la construction de l'aspect cardinal du nombre[2]. En particulier, les élèves doivent construire le nombre comme un outil de résolution de problèmes relatifs au dénombrement, à la comparaison ou la constitution de collections. Pour cela, un travail sur l'énumération (Briand *et al.*, 1999) est indispensable. Par ailleurs, à la fin de l'école maternelle, l'élève doit disposer de plusieurs procédures de dénombrement : le subitizing (perception globale de petites quantités), l'usage de constellations qui repose sur l'apprentissage d'associations (disposition d'objets - nombre), le comptage qui repose sur la connaissance de la comptine numérique et de procédures d'énumération, le calcul qui repose sur la connaissance de premiers faits numériques en lien avec les décompositions (si je perçois qu'il y a trois objets et encore deux, je sais qu'il y en cinq en tout) (Brissiaud, 2003, 2007). À cet égard, il faut noter que le travail des décompositions est indispensable à plusieurs titres. Les décompositions permettent de désigner des nombres de façon tout à fait pertinente dans des situations de l'école maternelle, en particulier pour comparer des collections (par exemple, repérer que dans l'une, il y a 3 et encore 2 et dans l'autre, 3 et encore 3, permet de conclure que la seconde est plus nombreuse). À l'école élémentaire, les décompositions sont aussi importantes pour le calcul (par exemple, dans l'addition posée « 35+27 », envisager « 5+7 » comme « 5 et encore 5 et encore 2 » permet de conclure qu'il y a « 2 unités et 3+2+1 dizaines »). Apprendre le comptage est tout aussi essentiel. En effet, cette procédure qui est toujours valable pour dénombrer une collection sera indispensable dans les classes ultérieures, notamment en lien avec le système de numération décimale de position. Par ailleurs, elle mobilise un principe essentiel de la numération : le fait que le nombre suivant de la comptine correspond à un de plus que le précédent. En variant les situations et les procédures mobilisées, à travers l'ensemble des situations qu'il rencontrera, l'élève construira un premier répertoire de faits numériques.

Concernant les aspects géométriques, la synthèse est plus difficile de notre point de vue. Les programmes insistent sur reconnaître, classer, nommer des formes ou des solides. Apprendre à distinguer, en situation fonctionnelle, un triangle d'un carré ou d'un disque est certainement important. Mais nous avons souvent observé en classe des travaux assez descriptifs (on compte le nombre de côtés et de « pics ») et définitoires. Or, apprendre ce qu'est un triangle, par exemple, peut aussi passer par l'expérimentation du fait qu'un triangle, contrairement au carré ou au rectangle, ne se déforme pas. En effet, si on construit un triangle avec des planchettes d'une certaine longueur vissées entre elles, les planchettes ne bougent pas et les angles sont préservés, ce qui n'est pas le cas avec un carré ou un rectangle fabriqué avec du matériel analogue. La maternelle peut aussi être l'occasion d'une première rencontre avec l'inégalité triangulaire : avec une très grande baguette et deux très petites on ne peut pas former un triangle !

---

[2] Pour cette raison, nous développons particulièrement cet aspect cardinal ici, au détriment de l'ordinal.

## 1.2. Installer un certain rapport aux mathématiques

L'école maternelle va inscrire l'élève dans un premier contrat didactique (Brousseau, 1998) en mathématiques et dans un premier rapport aux savoirs mathématiques (Charlot, 1997). Même si nous pensons que ce rapport peut évoluer dans la suite de la scolarité de l'élève, il semble important de lui prêter une attention particulière. Pour nous, il est essentiel de faire en sorte que l'élève n'envisage pas les mathématiques comme une discipline où une proposition est décrétée correcte ou pas « parce que l'enseignant l'a dit », ni où réussir consiste seulement à faire comme l'enseignant l'a dit ou montré. Cela implique, en particulier, que les situations proposées permettent, autant que possible, une validation au sens de la théorie des situations (Margolinas, 1993) : lorsqu'un élève propose une réponse, la situation renvoie des informations qui lui permettent de constater la validité ou la non validité de sa réponse. On retrouve là notre hypothèse d'apprentissage par adaptation. En ce sens, la validation s'oppose à l'évaluation qui est, elle, un avis de l'enseignant sur la production de l'élève.

Par ailleurs, fondamentalement, nous pensons que les savoirs mathématiques à l'école maternelle doivent être perçus comme des outils qui donnent un pouvoir d'agir sur les situations. Cela nous semble être l'une des conditions qui permet à un élève de s'engager avec intérêt et plaisir dans les situations mathématiques proposées à la maternelle mais aussi tout au long de sa scolarité. En effet, nous avons tous pu observer de jeunes élèves qui persévèrent quelquefois de longues minutes dans un problème, à côté d'un enseignant qui les encourage et régule leur activité sans pour autant réduire la difficulté, et se réjouissent largement lorsqu'ils parviennent à la solution. On attribue souvent cette persévérance et ce plaisir à l'aspect ludique de la situation proposée. C'est une part possible de l'explication. Mais on observe aussi souvent que ces élèves demandent à recommencer et sont tout aussi satisfaits de leur réussite. Il nous semble que si le plaisir de « jouer » ne s'émousse pas chez les jeunes élèves alors même que l'incertitude sur leur réussite diminue, c'est parce que, précisément, ils se réjouissent du pouvoir d'agir et de maîtriser la situation qu'ils sont en train d'acquérir.

Pour répondre à ces objectifs, les situations proposées aux élèves doivent remplir certaines conditions.

## 2. Le choix des situations

### 2.1. Pourquoi travailler sur fiche n'est pas satisfaisant ?

L'usage de fiches à l'école maternelle répond, pour partie, à un besoin de disposer de traces de l'activité des élèves lorsqu'ils sont en ateliers autonomes, traces utiles aux enseignants et disponibles pour les parents (Joigneaux, 2009). Les sociologues montrent que cet usage des fiches est source d'inégalités d'apprentissage, en particulier parce que les fiches impliquent un rapport à l'écrit et à son organisation qui est socialement discriminant (Joigneaux, 2015). En adoptant un regard didactique, nous allons montrer pourquoi un tel travail n'est pas satisfaisant pour les apprentissages des élèves malgré l'intérêt qu'il peut présenter en termes de traces. Pour cela, comparons deux situations de constitution d'une collection équipotente à une collection donnée qui sont quelquefois utilisées en maternelle avec l'objectif de permettre des apprentissages numériques. Nous effectuons pour cela une analyse de l'activité possible de l'élève. Ces deux activités peuvent être réalisées en atelier.

La première activité est une activité sur fiche (voir annexe 1) : une collection est constituée de ronds à l'intérieur d'un serpent, une autre de la représentation de jetons, on demande à l'élève de

« *colorier « autant » de jetons qu'il y a de ronds sur le serpent* ». Dans la seconde activité « *Juste assez* » (Thomas & Hersant, 2015) (voir annexe 2), les élèves disposent d'une boîte séparée en deux par un rabat qui permet de cacher l'une ou l'autre moitié de la boite. Des fiches avec des ronds dessinés dessus sont disposées en tas d'un côté de la boite, de l'autre côté de la boîte se trouvent des jetons. On demande à l'élève de « *prendre des jetons de façon qu'il y ait un jeton sur chaque rond, pas de jeton sans rond et pas de rond sans jeton* ».

### *Quelle activité mathématique des élèves dans chaque cas ?*

Dans le travail sur fiche, la consigne proposée mobilise le terme « autant ». Si ce terme n'est pas connu de l'élève, ce dernier ne peut rien faire, à part répondre au hasard et faire du coloriage. Il ne fait donc pas de mathématiques. En revanche, si l'élève connaît ce terme, il est fort probable que la solution soit immédiate pour lui. Cette activité sur fiche ne permet donc pas d'apprentissage. Tout au plus, elle peut servir pour évaluer certaines connaissances anciennes. En revanche, la consigne de la seconde activité n'utilise pas de mot savant, tous les élèves peuvent faire des essais, des erreurs et petit à petit comprendre que le nombre est l'outil pertinent pour répondre à la question.

Par ailleurs, le mot « autant » est souvent, pour les élèves, un indicateur qu'il faut compter. Dans l'activité 1, les élèves risquent donc de mobiliser le nombre plus par contrat didactique (Brousseau, 1998), pour répondre aux attentes de l'enseignant, que par nécessité. Dans l'activité 2, aucune indication n'est donnée ; cela doit amener l'élève à mobiliser de lui-même le nombre comme outil pour répondre au problème.

Enfin, même si comme nous l'avons indiqué, l'usage de mot « autant » dans la consigne engendre de façon quasi automatique un dénombrement chez les élèves, il faut noter que la présence de la première collection tout au long de l'activité peut amener à des solutions correctes sans que l'élève ne recoure au nombre. L'élève peut en effet très bien faire une correspondance terme à terme entre les deux collections. Cela n'est pas possible dans la seconde activité puisque la collection de départ est cachée au moment où l'élève choisit des jetons.

### *Fiche ou situation : qu'est-ce que ça change ?*

L'analyse précédente peut donner l'impression que la limite provient plus de la consigne proposée que de la modalité de travail sur fiche, puisque dans l'activité sur fiche, la compréhension de la consigne requiert la connaissance visée. Peut-on faire autrement ? Une consigne possible serait : « colorie les jetons pour qu'il y ait un jeton par rond, pas de rond sans jeton ». Avec cette consigne, puisque la première collection est toujours présente, les élèves risquent de faire une correspondance terme à terme sans mobiliser le cardinal de la collection. Ils peuvent donc réussir sans utiliser le nombre. Par ailleurs, on peut difficilement indiquer la condition associée « pas de jetons sans rond », il faudrait dire « pas de jeton colorié sans rond ». L'expression de ces deux conditions reste compliquée car les collections de jetons et de ronds ne sont pas déplaçables. À l'inverse, dans l'activité 2, la consigne est formulée sans recourir à l'idée de nombre. Cela est possible car on dispose d'une collection que l'on peut cacher et d'une collection d'objets déplaçables (les jetons), ce qui permet d'énoncer un but à l'activité (mettre un jeton sur chaque rond) ainsi qu'un critère de réussite[3] (pas de jeton sans rond, pas de rond sans jeton) sans évoquer les connaissances en jeu. L'activité de l'élève peut alors se définir à l'aide de ce but et de ce critère de réussite.

---

[3] Ici, le critère de réussite est un critère pour l'élève et non pour l'enseignant : si l'élève a posé un un jeton sur chaque rond et qu'il n'y a pas de rond sans jeton, alors il sait qu'il a réussi.

Par ailleurs, la fiche ne donne pas l'opportunité de recommencer si l'on a commis une erreur. En revanche, dans l'activité 2, l'élève peut recommencer autant de fois qu'il le souhaite avec un même fond ou un fond différent si l'enseignant l'y autorise. Or on sait qu'il est nécessaire de se confronter plusieurs fois à une situation pour apprendre.

Ensuite, dans la situation sur fiche, il n'y a pas de possibilité de validation au sens où nous l'avons défini précédemment. C'est l'enseignant qui indique si la production (collection de jetons coloriés) est correcte ou pas. Cette activité ne remplit donc pas notre hypothèse d'apprentissage par adaptation. Dans l'activité 2, l'élève sait, sans intervention du maître, que sa réponse convient ou pas, il n'a pas à l'accepter car le maître le dit.

On pourrait objecter que, dans l'activité 1, il est possible de proposer à l'élève un transparent pour vérifier sa proposition. Cela ne convient pas car il y a plusieurs façons de choisir 7 jetons parmi les 10 proposés. Quand bien même, on trouverait une possibilité avec un transparent, cette vérification n'aurait pas, à notre avis, la même force que celle qui se réalise dans l'activité 2. En effet, c'est l'enseignant qui a produit le transparent et alors la réponse n'est pas juste ou fausse en lien avec une nécessité matérielle et mathématique mais en lien avec l'acceptation de ce que dit le maître.

L'activité 2 remplit les principales conditions que Douady (1984) pose pour une situation-problème et que nous allons détailler ci-dessous.

## 2.2. Quelles caractéristiques pour les situations proposées aux élèves ?

### *Proposer des situations-problèmes*

Douady propose les caractéristiques suivantes pour une situation-problème1[4] (Douady, 1984) :

1. l'élève doit pouvoir s'engager dans la résolution du problème. Il peut envisager ce qu'est une réponse possible du problème ;
2. les connaissances dont il dispose (au début de l'activité) sont insuffisantes pour qu'il résolve immédiatement le problème ;
3. la situation problème doit permettre à l'élève de décider si une solution trouvée est convenable ou pas ;
4. la connaissance que l'on désire voir acquérir par l'élève doit être l'outil le plus adapté pour la résolution du problème au niveau de l'élève.

Commentons rapidement ces critères.

Le premier nous indique qu'il ne faut pas que l'élève reste « sec » face au problème posé, sinon il n'apprendra pas et il risque de se décourager. On peut aussi se trouver dans une situation où l'élève cherche à deviner les attentes et intentions de l'enseignant ou des indices de solution dans ce que dit l'enseignant. Et, dans ce cas, l'élève ne fait pas de mathématiques ; il risque de s'inscrire dans une recherche de conformité à ce que l'enseignant attend. Pour remplir le premier critère énoncé par Douady, le choix du problème et du matériel mis à disposition des élèves pour y répondre est essentiel. Mais comme nous le verrons dans le paragraphe suivant, la façon dont l'enseignant introduit la situation a aussi toute son importance.

Nous ne nous attarderons pas sur le second critère qui nous semble reprendre simplement l'idée de problème. On peut noter toutefois que dans le travail sur fiche, le plus souvent, ce critère essentiel n'est pas rempli, comme nous l'avons montré dans le paragraphe précédent.

---

[4] Douady ajoute une cinquième condition sur les changements de cadre : « *Le problème peut se formuler dans plusieurs cadres et entre lesquels on peut établir des correspondances* ». Cette condition nous semble difficile à remplir à l'école maternelle.

Le troisième critère renvoie à la question de la validation. Là encore, le choix du problème et du matériel va intervenir de façon cruciale pour le respect de ce critère. Mais, les régulations de l'enseignant vont avoir leur importance. D'abord, parce qu'il peut court-circuiter cette validation s'il intervient trop rapidement pour indiquer à l'élève si ce qu'il propose est correct ou pas. Ensuite, parce que, en maternelle, les élèves peuvent avoir besoin d'une aide pour estimer, au début, la validité de leur réponse à un problème. L'enjeu est alors d'intervenir sans court-circuiter la validation, tout en permettant de faire comprendre que la réponse apportée n'est pas satisfaisante. De ce point de vue, avoir défini des critères de réussite est très utile. Par exemple, dans une situation classique de dénombrement, l'enseignant pourra dire : « Va chercher des jetons pour qu'il y ait un jeton sur chaque rond. Tu auras réussi s'il n'y a pas jeton sans rond et pas de rond sans jeton. ». Lorsqu'il ramène ses jetons, si l'élève ne conclut pas lui-même sur la validité de son résultat, l'enseignant pourra lui rappeler le critère de réussite pour qu'il s'engage dans une validation.

Le quatrième critère peut paraître une évidence, sinon effectivement l'élève passe à côté des apprentissages visés. Les données du problème, le matériel que l'on met à disposition sont des variables qui influencent l'activité des élèves. Quelquefois, des « détails » peuvent avoir des conséquences importantes sur l'activité des élèves. Ainsi, par exemple, si on veut amener les élèves à utiliser le nombre pour constituer une collection équipotente à une collection donnée, il faut que la collection de référence ne soit plus visible lorsque l'élève constitue la nouvelle collection. Sans quoi l'élève peut simplement effectuer une correspondance terme à terme. Il faut donc soit éloigner la collection de référence des jetons à utiliser pour constituer la collection équipotente — c'est la situation fondamentale du dénombrement pour Brousseau (1997) — soit cacher cette collection comme dans la situation 2. Avant de proposer la situation, se demander « que peut faire l'élève ? » et y répondre en envisageant des possibles est souvent utile pour s'assurer que la situation proposée contraint bien l'élève à mobiliser le savoir visé.

Ces critères peuvent servir de guide pour se faire une idée sur l'intérêt didactique d'une situation proposée dans une ressource. Ils peuvent aussi constituer des points d'appui pour modifier une situation proposée dans une ressource. Reprenons par exemple l'activité 1. On peut facilement la transformer pour obtenir une situation-problème. Il suffit d'enlever la partie basse de la feuille où figurent les jetons à colorier, de mettre à disposition des vrais jetons, d'éloigner les deux collections de façon qu'elles ne soient pas visibles simultanément et de modifier la consigne.

*Permettre à l'élève d'entrer dans le problème*

Revenons sur la première caractéristique d'une situation-problème pour préciser le rôle que l'enseignant peut jouer afin que l'élève s'engage dans la situation et envisage une réponse possible au problème. Comparons pour cela deux réalisations d'une même situation, *« Le bon pot »* (Thomas & Hersant, 2015) par deux enseignants débutants dans deux classes de PS.

La situation est la suivante. On dispose de 4 pots retournés et d'une petite peluche. Sur chaque pot est collé un couvercle dans lequel sont disposés des cubes en très petite quantité. L'enseignant :

1. cache la peluche sous un des pots et annonce : « Je cache la peluche sous le pot à deux cubes »,
2. après avoir demandé aux élèves de fermer les yeux, change les pots de place en les

faisant glisser, sans les soulever.

La peluche reste ainsi sous le même pot mais le mouvement fait bouger la position des cubes sur les couvercles. Les élèves doivent désigner du doigt le pot sous lequel est cachée la peluche.

Dans les deux classes, la situation est proposée à des élèves de niveau hétérogène en mathématiques. Dans la première classe, tous les élèves sont investis dans la situation, ils ne cherchent pas à soulever les pots pour retrouver la peluche, sollicitent l'enseignante pour « jouer » et prennent, à l'évidence, plaisir à « parier » sur le bon pot. Dans l'autre classe, les élèves s'engagent peu dans le pari et n'ont de cesse d'essayer de soulever tous les pots pour retrouver la peluche.

Comment expliquer cette différence d'attitude ? Plusieurs variables peuvent intervenir. Mais, les deux enseignantes utilisent un matériel semblable à celui proposé dans la ressource et formulent leur question aux élèves de façon analogue à ce qui est proposé dans la ressource. On peut donc penser que la différence ne se situe pas là. En fait, les deux enseignantes introduisent de façons très différentes le problème. Dans la première classe, l'enseignante propose une première phase où elle mélange les pots sous les yeux des élèves puis une seconde phase où elle demande aux élèves de fermer les yeux. Dans la seconde classe, l'enseignante propose de fermer les yeux dès la phase 1.

Nous pensons que, dans la première classe, la phase 1 permet aux élèves de comprendre :
1. qu'il n'y a aucune raison pour que la peluche soit au même endroit après le mélange ;
2. que l'enseignant fait seulement glisser les pots, il ne les soulève pas et ne change pas la peluche de pot.

Ils peuvent alors imaginer qu'il ne s'agit pas d'une question de place des pots, ni de hasard et qu'il faut trouver un autre critère pour résoudre le problème. Cela aide les élèves à envisager une réponse possible dans la mesure où cela leur permet d'écarter certains possibles. Pour le dire autrement, cela leur permet de construire le problème (Fabre & Orange, 1997 ; Hersant, 2022) et de s'en faire une représentation (Julo, 1995). Dans la seconde classe, les élèves ne nous semblent pas avoir cette opportunité et sont alors réduits à essayer, au hasard, pour retrouver la peluche.

Cet exemple illustre le rôle de l'enseignant dans la dévolution d'une situation-problème. Le problème choisi, le matériel mis à disposition sont importants mais ils ne suffisent pas à définir une situation-problème. La façon dont l'enseignant présente le problème a un rôle déterminant pour permettre à l'élève d'en imaginer une solution. Elle intervient aussi de façon importante pour installer une situation sur laquelle l'élève sent qu'il peut agir et qu'il a les moyens de « gagner » à l'intérieur d'un cadre clairement défini.

*Proposer plusieurs fois une situation et formuler les savoirs*

À la maternelle, proposer aux élèves des situations-problèmes n'est une option intéressante que si on leur donne la possibilité de s'y confronter plusieurs fois. En effet, confronté à une situation-problème, un élève est en situation de construire des connaissances. Mais on ne trouve pas toujours du premier coup la solution, on peut trouver une solution qui ne fonctionne que partiellement et qu'il va donc falloir affiner. La construction des connaissances pérennes demande du temps. Il est donc souhaitable que l'élève « joue » à plusieurs reprises au cours d'une séance et qu'une situation soit proposée plusieurs fois, dans une forme analogue mais aussi avec des « habillages » qui évoluent.

De plus, nous l'avons indiqué dès le début de cet article : l'élève n'apprend pas seulement en se confrontant à des situations, même lorsqu'elles correspondent à des situations-problèmes, il apprend aussi parce que l'enseignant organise un processus d'institutionnalisation (Brousseau, 1998) qui va permettre de donner aux connaissances produites en situation par les élèves un statut de savoir, partagé dans la classe, et au-delà, exigible dans la scolarité future. Mener à bien ce processus d'institutionnalisation nécessite d'identifier les savoirs et d'envisager des formulations de ces savoirs. Il est souvent, à tort, envisagé seulement comme une phase : celle dans laquelle l'enseignant explique ce qu'il y a à retenir. Or, indiquer à la fin de la situation ce qu'il faut en retenir ne peut être satisfaisant. En effet, d'abord cela ne permet pas, le plus souvent, de tisser les relations entre les connaissances mises en œuvre par les élèves en situation (leurs découvertes) et ce qu'on souhaite que les élèves retiennent de la situation, ce qui est essentiel pour que les savoirs soient fonctionnels. Proposer des formulations tout au long de la situation permet entre autres de dire plusieurs fois de façons différentes, plus ou moins décontextualisées et plus ou moins comme des stratégies, les savoirs en jeu dans la situation. Il s'agit évidemment de proposer des formulations sans dévoiler trop tôt la procédure qui permet de réussir et ainsi priver la plupart des élèves de la construction de cette procédure. Mais il est aussi important que tous les élèves accèdent, à terme, aux procédures gagnantes. Il y a donc un subtil équilibre à trouver entre adaptation et acculturation, qui renvoie à une dialectique de la dévolution et de l'institutionnalisation.

*Tout ne s'enseigne pas à l'aide de situations-problèmes*

Les situations-problèmes offrent de bonnes conditions pour permettre des apprentissages mathématiques. Cependant, toutes les notions mathématiques ne se prêtent pas à l'élaboration d'une situation-problème. En particulier, tous les savoirs associés à des conventions ne peuvent pas, en tant que tels, faire l'objet d'un apprentissage par résolution de problème. À l'école maternelle, l'apprentissage, crucial, de la comptine numérique relève de cette catégorie. En effet, il s'agit d'apprendre une suite ordonnée et conventionnelle de mots. On ne peut pas demander aux élèves de construire ces mots et cet ordre. Il donc faut trouver des moyens ludiques de faire répéter les élèves tout en leur faisant comprendre la nécessité de cet apprentissage pour les situations de dénombrement. Il en est de même pour l'apprentissage de l'écriture chiffrée des premiers nombres. L'élève doit apprendre à reconnaître que le signe « 4 » correspond à la quantité quatre mais on a intérêt à lui faire ressentir la nécessité et l'intérêt de cet apprentissage dans des situations où l'échange d'informations sur les quantités ne se fait que par écrit et à lui montrer la puissance de cette connaissance pour garder sur un temps long la mémoire d'une quantité.

## 3. Identifier les fonctions du matériel et des phases d'une situation pour réguler l'activité des élèves

Dans la partie **1.** de cet article, à propos de la comparaison d'un travail sur fiche et d'un travail en situation, nous avons montré que le matériel mis à disposition des élèves dans *« Juste assez »* conférait à cette situation des propriétés difficiles à retrouver dans une activité sur fiche. Il ne faudrait pas raccourcir cette idée en pensant que la différence entre les deux activités tient au fait que, dans un cas, l'élève manipule et pas dans l'autre. Il est aussi essentiel en effet de percevoir la fonction du matériel dans une situation.

## 3.1. Les fonctions didactiques du matériel…

Pour illustrer notre propos, nous allons nous appuyer sur la comparaison de deux situations, *« Juste assez »* (situation 2 dans ce texte) et *« Les deux cartes »* (Thomas & Hersant, 2015).

Cette dernière situation vise à amener les élèves à utiliser le dénombrement pour comparer et apparier deux collections non déplaçables. Les élèves ont à disposition des cartes sur lesquelles sont représentés des ronds (pas forcément organisés en constellation, ni en ligne), des jetons dans une boîte qui ferme et des épingles à linge. Le problème consiste à choisir une carte puis à l'apparier avec une autre qui a autant de ronds. La consigne est la suivante : « Il faut choisir des cartes qui ont autant[5] de ronds ». Elle permet de poser le problème sans donner d'indication sur la ou les procédures pertinentes pour le résoudre. Pour réussir, l'élève doit comparer les collections, ce qui implique de dénombrer la première collection, de garder en mémoire son cardinal et de comparer ce nombre avec le cardinal des collections représentées sur d'autres cartes. Selon la disposition des ronds sur les cartes et selon le nombre de ronds, l'élève peut utiliser différentes procédures pour dénombrer les collections. Lorsqu'il pense avoir trouvé une carte de même cardinal que la première carte choisie, il vérifie à l'aide de jetons : il pioche dans la boîte autant de jetons que nécessaire pour recouvrir tous les ronds d'une des cartes puis ferme la boîte et dispose les jetons sortis sur les ronds de la seconde carte. S'il n'y pas de jeton sans rond, ni de rond sans jeton, l'élève a bien reconnu que les deux cartes ont le même cardinal, il les apparie alors avec une épingle à linge. Ce critère de réussite « pas de jeton sans rond et pas de rond sans jeton » est indiqué aux élèves au moment de la présentation du problème.

Le matériel utilisé dans cette situation est assez semblable à celui de *« Juste assez »* : les cartes avec des ronds dessinés et des jetons sont communs. Mais dans ces deux situations, le matériel n'a pas la même fonction. Dans *« Juste assez »*, les cartes correspondent aux collections de départ et les jetons permettent la constitution des collections équipotentes. Ce matériel permet à l'élève de faire des tentatives et de s'engager dans le problème. Par ailleurs, ce même matériel a aussi un rôle de validation : lorsque l'élève dispose ses jetons sur les ronds de la carte, il sait s'il a réussi ou pas. Ce sont ces informations, et les interventions de l'enseignant dans cette phase de validation, qui vont lui permettre, petit à petit, de construire les nécessités associées au problème (il y a un rond sans jeton, il aurait fallu que je prenne un jeton de plus ; si j'ai 5 ronds, il faut 5 jetons). Dans « Les deux cartes », les cartes à apparier correspondent aux collections à associer. Les jetons interviennent seulement dans la validation. Cette validation renvoie à l'élève des informations qui lui permettent de construire le problème.

La distinction de ces fonctions du matériel est essentielle pour réguler l'activité des élèves.

## 3.2. … Et la régulation de l'activité des élèves

Dans la situation *« Les deux cartes »*, considérons le cas d'un élève qui connaît parfaitement sa comptine dans le champ numérique concerné et dénombre sans difficulté les collections de ronds. Nous considérons que, pour cet élève, l'enjeu est de saisir que le dénombrement des collections est la connaissance la plus adaptée pour résoudre le problème. Après quelques parties, il est fort probable que l'élève court-circuitera la phase de validation avec les jetons et s'affairera à constituer le plus de paires de cartes possibles. Il n'y a alors pas de raison de l'obliger à utiliser les jetons pour vérifier que les deux cartes ont bien le même nombre de ronds.

À l'inverse, pour un élève qui peine à trouver comment associer deux cartes, l'utilisation de collections de jetons est essentielle pour qu'il identifie que sa production ne remplit pas les critères de réussite fixés. Il est donc important qu'il continue d'utiliser les jetons pour vérifier son résultat. Le rôle de l'enseignant peut alors consister à mettre en évidence que la réponse de l'élève ne remplit pas les critères de réussite posés et à pointer les différences (par exemple : est-ce que tu as un jeton sur chaque rond ?). En revanche faire compter à l'élève le nombre de jetons ou introduire des informations numériques à propos des jetons (par exemple : il te manque deux jetons) ne présente aucun intérêt (Hersant, 2020). Au contraire, cela le détourne du problème qu'on souhaite qu'il construise, voire l'incite à mobiliser le nombre sans ressentir la nécessité de son usage. En revanche, dans *« Juste assez »*, au moment de la validation, des informations comme « il devrait y avoir un jeton sur chaque rond mais il y a deux jetons sans ronds » apparaissent à un moment incontournables pour aider l'élève à se représenter le problème.

---

[5] Cette situation est proposée que lorsque le terme « autant » est compris par les élèves.

## Conclusion

À travers ces exemples, nous avons cherché à mettre en évidence des conditions qui permettent une activité source d'apprentissage chez les élèves. Ces conditions sont principalement organisées autour des concepts-clés de situation-problème (Douady, 1984) et de validation (Brousseau, 1998) ; elles renvoient pour la majorité au travail qu'effectue l'enseignant en amont de la réalisation de la situation en classe : identification des savoirs visés, conception et choix d'une situation dont les caractéristiques contraignent à mobiliser le savoir visé et permettent la validation, identification des fonctions didactiques du matériel, formulation d'une consigne qui précise un but et des critères de réussite, anticipation de la façon dont l'activité va être présentée aux élèves et des interventions auprès des élèves au cours de l'activité. Toutes ces conditions demandent à l'enseignant d'anticiper avec beaucoup de précision le déroulement de la situation. Ainsi, même si dans cet article nous avons soutenu notre propos à partir de l'analyse de l'activité d'élèves en situation, on le voit bien, l'activité de l'élève dépend fortement du travail de l'enseignant en amont de la séance.

Nous l'avons souligné à plusieurs endroits dans ce texte, la possibilité d'une activité mathématique des élèves à la maternelle réside pour partie dans la mise en place de conditions qui permettent à l'élève d'exercer un pouvoir d'agir sur les situations. Installer ces conditions, c'est d'abord permettre des apprentissages mathématiques. C'est aussi une façon de montrer aux jeunes élèves à quoi servent les mathématiques et une façon de donner à tous, filles ou garçons, quelle que soit la classe sociale dont ils·elles sont issu·e·s, le goût des mathématiques. Surtout, fondamentalement, c'est installer, dès la maternelle, un rapport aux mathématiques qui n'est pas celui d'une discipline où on « bloque », où on ne comprend pas ce qu'il « faut faire » ou ce que la « maîtresse attend » mais celui d'une discipline qui, au contraire, permet une émancipation par les savoirs. En effet, en mathématiques, il ne s'agit pas de résoudre des problèmes pour faire plaisir à l'enseignant, la plupart du temps ce sont les nécessités de la situation qui déterminent la validité d'une réponse.

## Références bibliographiques

Valentin, D. (2004). *Découvrir le monde avec les mathématiques : Situations pour la petite et la moyenne sections*. Hatier.

MEN (2021). BOENJS n°25 du 24 juin 2021

**Annexe 1
Fiche**

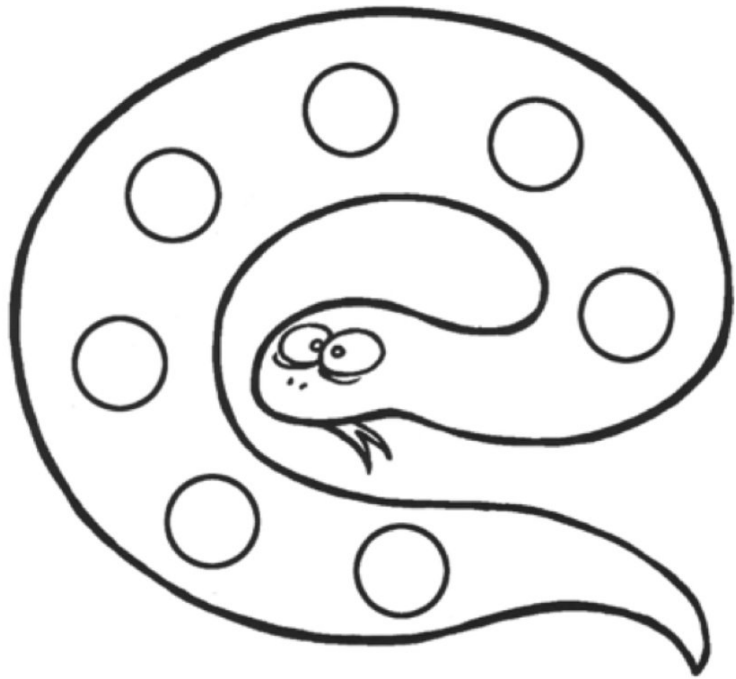

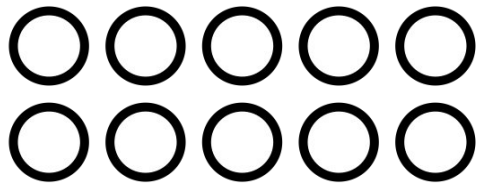

# Annexe 2
# Situation « juste assez »

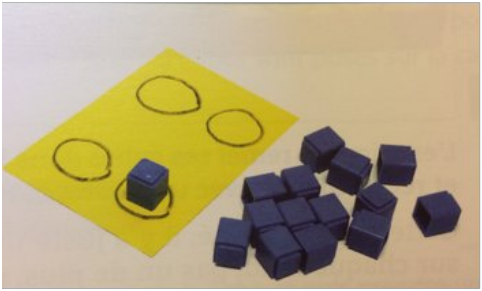
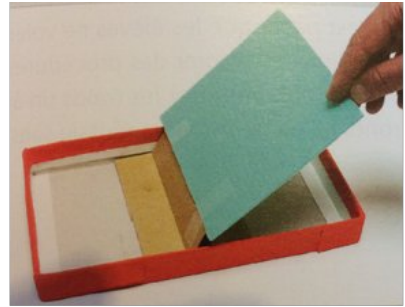
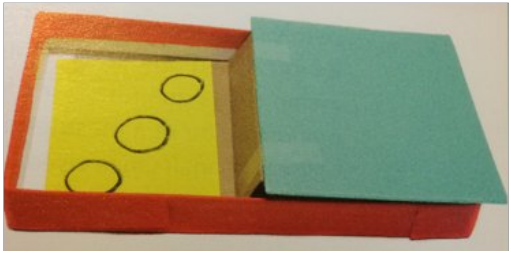
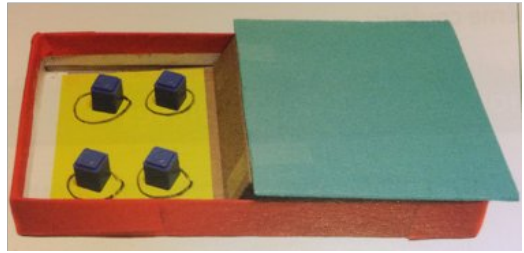